\documentclass[a4paper,fleqn,10pt]{amsart}

\usepackage{natbib}
\usepackage{mathscinet}
\usepackage{mathrsfs}
\usepackage{hyperref}
\hypersetup{
    pdftitle={On the infimum attained by a reflected L\'evy process},
    pdfsubject={Probability Theory},
    pdfauthor={Krzysztof Debicki and Kamil Kosinski and Michel Mandjes},
    pdfdisplaydoctitle=true,
}

\newcommand{\ee}{\mathbb E}

\newcommand{\nn}{\mathbb N}
\newcommand{\pp}{\mathbb P}
\newcommand{\rr}{\mathbb R}

\newcommand{\prob}[1]{\pp\left( #1 \right)}
\newcommand{\as}[1]{\quad\text{as}\quad #1\toi}

\newcommand{\toi}{\to\infty}
\newcommand{\de}{\,{=}_{\rm d}\,}
\newcommand{\clasS}{\mathscr S}
\newcommand{\D}{\mathrm d}

\newtheorem{thm}{Theorem}

\newtheorem{cor}{Corollary}
\newtheorem{prop}{Proposition}
\newtheorem{assumption}{Assumption}
\theoremstyle{remark}
\newtheorem{rem}{Remark}

\newcommand{\vb}{\vspace{1.5mm}}

\begin{document}
\bibliographystyle{plain}
\setcitestyle{numbers}

\title
{On the infimum attained by a reflected L\'evy process}

\author{K.\ D\polhk{e}bicki}
\address{Instytut Matematyczny, University of
Wroc\l aw, pl.\ Grunwaldzki 2/4, 50-384 Wroc\l aw, Poland.}
\email{Krzysztof.Debicki@math.uni.wroc.pl}
\thanks{KD was supported by
MNiSW Grant N N201 394137 (2009-2011) and by a travel grant from
NWO (Mathematics Cluster STAR)}

\author{K.M.\ Kosi\'nski}
\address{Korteweg-de Vries Institute for Mathematics,
University of Amsterdam, the Netherlands; E{\sc urandom},
Eindhoven University of Technology}
\email{K.M.Kosinski@uva.nl}
\thanks{KK was supported by NWO grant 613.000.701.}

\author{M.\ Mandjes}
\address{Korteweg-de Vries Institute for Mathematics,
University of Amsterdam, the Netherlands; E{\sc urandom},
Eindhoven University of Technology, the Netherlands; CWI,
Amsterdam, the Netherlands}
\email{M.R.H.Mandjes@uva.nl}
\thanks{KD and MM thank the Isaac Newton Institute, Cambridge, for hospitality.}


\date{\today}

\begin{abstract}
This paper considers a L\'evy-driven queue (i.e., a L\'evy process
reflected at 0), and focuses on the distribution of
$M(t)$, that is, the minimal value attained in an interval of
length $t$ (where it is assumed that the queue is in stationarity
at the beginning of the interval). The first contribution is an
explicit characterization of this distribution, in terms of
Laplace transforms, for spectrally one-sided L\'evy processes
(i.e., either only positive jumps or only negative jumps). The
second contribution concerns the asymptotics of $\prob{M(T_u)> u}$
(for different classes of functions $T_u$ and $u$ large); here we have to distinguish
between heavy-tailed and light-tailed scenarios.
\end{abstract}

\maketitle

\section{Introduction}
\label{sec:intro}
The class of processes with stationary and independent increments,
known as {\it L\'evy processes}, form a key object in applied
probability. A substantial body of literature is devoted to L\'evy
processes that are reflected at 0, sometimes also referred to as
{\it L\'evy-driven queues}, and are regarded as a valuable
generalization of the classical M/G/1 queues; also, the important
special case of reflected Brownian motion is covered.

These reflected L\'evy processes are defined as follows. Let
$X\equiv\{X(t):t\in\rr\}$ be a L\'evy process with (without loss of generality)
zero drift: $\ee X(1) = 0$ and $X(0)=0$. Then define the queueing process
(or: workload process, storage process) $Q\equiv\{Q(t):t\ge 0\}$ through (for $c>0$)
\[Q(t):= \sup_{s\le t} \left(X(t) - X(s) -c (t-s)\right),\]
where it is assumed that the workload is in {\it equilibrium} (stationarity) at time 0, i.e., $Q(0)\de Q_e$.

We refer to this process $Q$ as the reflection of the
L\'evy process $Y=\{Y(t):t\in\rr\}$ at 0, where $Y(t):=X(t)-ct$. In the
sequel we normalize time such that $c=1$.

\vb

When considering the steady state $Q_e$ of the reflected process
introduced above, the literature can be roughly divided into two
categories. (A)~In the first place there are results on the full
distribution of $Q_e$, in terms of the corresponding Laplace
transform. Particularly for the case of {\it one-sided jumps},
these transforms are fairly explicit. 
If $X$ is such that it has only positive jumps and is not a {\it subordinator} (an increasing process), $X\in\clasS^+$ (which is often referred to as the {\it
spectrally positive case}), then a generalization of the classical
Pollaczek--Khintchine formula was derived \cite{Zolotarev64}, while
in the case of only negative jumps, $X\in\clasS^-$ ({\it spectrally negative}),
$Q_e$ was seen to be exponentially distributed. In the L\'evy
processes literature \cite{Bertoin96,Kyprianou06}, this type of
results can be found under the denominator {\it fluctuation
theory}. We recall that there are powerful tools available for
numerical inversion of Laplace transforms \cite{Abate95, Iseger06}. (B)~In
the second place there are results that describe the {\it
asymptotics} of $\prob{Q_e>u}$ for $u$ large. Then one has to
distinguish between results in which the upper tail of the L\'evy
increments is light on the one hand, sometimes referred to as the
Cram\'er case, and results that correspond to the heavy-tailed
regime on the other hand; see for instance \cite{Dieker06} and
references therein.

In the present short communication, we consider a related problem:
we analyze how long the process consecutively spends above a given
level. More formally, we consider the distribution of
$M(t):=\inf_{s\in[0,t]} Q(s)$, i.e., the minimum value attained by
the workload process in a window of length $t$, where it is
assumed that the queue is in stationarity at the beginning of the
interval. This problem has various applications: one could for
instance think of the analysis of persistent overload in an
element of a communication network or a node in a supply chain;
see e.g.\ \cite{Pacheco08}. A related study on the situation of
infinitely-divisible self-similar input is \cite{Albin04}.

Our results correspond to both branches (A) and (B) mentioned
above: in \autoref{sec:transforms} we present results on the
Laplace transform of $M(t)$, relying on known results for L\'evy
fluctuation theory; we also consider the special case of Brownian
motion. \autoref{sec:asymptotics} identifies the asymptotics of
$\prob{M(T_u)>u}$ for different classes of functions $T_u$ and $u$
large; as expected, we need to distinguish between heavy-tailed
and light-tailed input.

Recall that
\[
Y(t)=X(t)-t, \quad Q(t)=\sup_{s\le t}\left(Y(t)-Y(s)\right),
\quad M(t)=\inf_{s\in[0,t]} Q(s).
\]
For a stochastic process $Z$ we will write:
\[
\underline{Z}(t):=\inf_{s\in[0,t]} Z(s), \quad \overline{Z}(t):=\sup_{s\in[0,t]} Z(s).
\]
It is well known that the process $Q$ admits the following representation
(see \cite[p. 375]{Robert03}):
\[
Q(t)=Q(0)+Y(t)+\max\left(0,-Q(0)-\underline{Y}(t)\right)
\]
so that, for any $u>0$, 
\[
\prob{M(t)>u}=\prob{Q(0)+\underline{Y}(t)>u}.
\]
Notice that due to the independent increments property of $X$, the random variables $Q(0)$ and $\underline{Y}(t)$
are independent, and hence $Q(0)+\underline{Y}(t)\de Q_e + \underline{Y}(t)$.

\section{Transforms for the spectrally one-sided case}
\label{sec:transforms}
\newcommand{\e}[1]{\boldsymbol{e}_{#1}}
Let $\e q$ denote a generic, exponentially distributed random variable with parameter $q>0$, that is independent from the process $X$.
In this section we evaluate the double transform, with $x,q> 0$,
\begin{equation}
\label{eq:L}
{\mathscr L}(x,q):= \ee e^{-xM(\e q)}=1 - x{\mathscr K}(x,q),
\end{equation}
where
\[
{\mathscr K}(x,q):=\int_0^\infty e^{-xu}\prob{M(\e q)>u}\D u
=\int_0^\infty e^{-xu}\prob{Q_e+\underline{Y}(\e q)>u}\D u
\]
and we used the fact, that for a non-negative random variable $K$ and any $x>0$, $\ee e^{-x K}=1-x\int_0^\infty e^{-x u}\prob{K>u}\D u$. As indicated in the introduction, we compute $\mathscr L$ for L\'evy processes with one-sided jumps. 

It is worth noticing that the double transform $\mathscr L$ uniquely determines the distribution of $M(t)$.
Determining the probability distributions from the double transform requires Laplace inversion. 
It is noted that recently, substantial progress has been made with respect to this type of inversion techniques. 
Besides the `classical' reference \cite{Abate95}, we wish to draw attention to significant recent progress in \cite{Iseger06}; 
the latter reference specifically addresses the multidimensional transforms, 
and also provides a fairly complete literature overview.
\subsection{Spectrally negative case}
\label{subsec:SNC}
Let $\psi(x)=\log\ee e^{x Y(1)}$ denote the Laplace exponent of $Y$ (which in this case is defined for every $x\ge 0$) and $\Phi(q)=\sup\{\theta\ge0:\psi(\theta)=q\}$
be the right inverse of $\psi$. It is well known, that in this case $Q_e$ is exponentially distributed with parameter $\Phi(0)>0$. Furthermore, the Wiener-Hopf factorisation, see \cite[Theorem 6.16]{Kyprianou06}, yields, for any $x,q>0$,
$
\ee e^{x\underline{Y}(\e q)}=\hat\kappa(q,0)/\hat\kappa(q,x),
$
where $\hat\kappa(x,y)=(x-\psi(y))/(\Phi(x)-y)$. That is,
\[
\ee e^{x\underline{Y}(\e q)}=\frac{q}{\Phi(q)}\frac{x-\Phi(q)}{\psi(x)-q}.
\]
\begin{thm}
\label{thm:SN}
For a spectrally negative process $X$,
\[
{\mathscr L}(x,q) =\frac{\Phi(0)}{x+\Phi(0)}\frac{\Phi(q)+x}{\Phi(q)},\quad x,q>0.
\]
\end{thm}
\begin{rem}
In the spectrally negative case, it is known that $\overline{Y}(\e q)$ is exponentially distributed with parameter $\Phi(q)$. Therefore, \autoref{thm:SN} can be reformulated as
\[
\ee \exp(-x M(\e q)) =\frac{\ee \exp(-xQ_e)}{\ee \exp(-x\overline Y(\e q))}.
\]
\end{rem}
\begin{proof}
Due to the independence of $Q_e$ and $\underline{Y}(\e q)$, and using the fact that $Q_e$ is exponentially distributed,
\begin{align*}
{\mathscr K}(x,q)&=
\int_0^\infty e^{-xu}
\int_{(-\infty,0]} \prob{Q_e>u-z} \D \prob{\underline{Y}(\e q)\le
z}\D u\\
&=
\int_0^\infty e^{-xu}
\int_{(-\infty,0]} e^{-\Phi(0)(u-z)} \D \prob{\underline{Y}(\e q)\le
z}\D u\\
&=
\int_0^\infty e^{-u(x+\Phi(0))}\D u
\int_{(-\infty,0]} e^{\Phi(0)z} \D \prob{\underline{Y}(\e q)\le
z}\\
&=
\frac{1}{x+\Phi(0)}\ee e^{\Phi(0)\underline{Y}(\e q)}
=
\frac{1}{x+\Phi(0)}\frac{\Phi(q)-\Phi(0)}{\Phi(q)}
\end{align*}
and the claim follows from \eqref{eq:L}.
\end{proof}

\subsection{Spectrally positive case}
Let $\hat\psi(x)=\log\ee e^{x \hat Y(1)}$
be the Laplace exponent of $\hat Y=-Y$ and $\hat\Phi(q)=\sup\{\theta\ge 0:\hat\psi(\theta)=q\}$ be the right inverse of $\hat\psi$. In this case, the well known Pollaczek--Khintchine formula gives, for $x>0$,
\[
\ee e^{-x Q_e}=\frac{\hat\psi'(0)x}{\hat\psi(x)}.
\]
Furthermore, note that $\underline{Y}(\e q)=-\overline{\hat Y}(\e q)$. It is known that
$\overline{\hat Y}(\e q)$ is exponentially distributed with parameter $\hat\Phi(q)$. 
\begin{thm}
\label{thm:SP}
For a spectrally positive L\'evy process $X$,
\[
{\mathscr L}(x,q) =
\frac{\hat\psi'(0)x}{\hat\psi(x)}\frac{\hat\Phi(q)}{q}\frac{q-\hat\psi(x)}{\hat\Phi(q)-x},\quad x,q>0,
\]
where the right hand side is understood in the asymptotic sense when $x=\hat\Phi(q)$, that is $\hat\psi'(0)\hat\Phi^2(q)\hat\psi'(\hat\Phi(q))/q^2$.
\end{thm}
\begin{rem}
By \citep[Theorem 4.8]{Kyprianou06} (understood in the asymptotic sense as well
when $x=\hat\Phi(q)$),
\[
\ee e^{-x(\overline{\hat Y}(\e q)-\hat Y(\e q))}=
\ee e^{x\underline{\hat Y}(\e q)}=
\ee e^{-x\overline{Y}(\e q)}=\frac{q}{\hat\Phi(q)}\frac{x-\hat\Phi(q)}{\hat\psi(x)-q}.
\]
Therefore, \autoref{thm:SP} can be reformulated as
\[
\ee \exp(-x M(\e q)) =\frac{\ee \exp(-xQ_e)}{\ee \exp(-x\overline Y(\e q))}.
\]
\end{rem}
\begin{proof}
First, let us assume that $x\ne\hat\Phi(q)$.
Due to the independence of $Q_e$ and $\underline{Y}(\e q)$, and
using the fact that $\overline{\hat Y}(\e q)$ is exponentially distributed,
\begin{align*}
{\mathscr K}(x,q)&=
\int_0^\infty e^{-xu}
\int_{[u,\infty)} \prob{\underline{Y}(\e q)>u-z} \D \prob{Q_e\le
z}\D u\\
&=\int_{[0,\infty)}\int_0^z e^{-xu}
 \prob{\underline{Y}(\e q)>u-z} \D u \D \prob{Q_e\le z}\\
&= 
\int_{[0,\infty)}e^{-xz}\int_0^z e^{wx}
 \prob{\overline{\hat Y}(\e q)<w} \D w \D \prob{Q_e\le z}\\
&=
\int_{[0,\infty)}e^{-xz}\int_0^z e^{wx}
\left(1-e^{-w\hat\Phi(q)}\right) \D w \D \prob{Q_e\le z}\\
&=
\frac{1}{x}\int_{[0,\infty)}(1-e^{-xz})\D \prob{Q_e\le z}
-
\frac{1}{x-\hat\Phi(q)}\int_{[0,\infty)}(e^{-\hat\Phi(q)z}-e^{-xz})\D \prob{Q_e\le z}\\
&=
\frac{1}{x}(1-\ee e^{-xQ_e})
-
\frac{1}{x-\hat\Phi(q)}(\ee e^{-\hat\Phi(q)Q_e}-\ee e^{-x Q_e})\\
&=
\frac{1}{x}\left(1-\frac{\hat\psi'(0)x}{\hat\psi(x)}\right)-
\frac{\hat\psi'(0)}{(x-\hat\Phi(q))}\left(\frac{\hat\Phi(q)}{q}-\frac{x}{\hat\psi(x)}\right).
\end{align*}

If $x=\hat\Phi(q)$, then the same computations give
\[
\mathscr K(\Phi(q),q)=
\frac{1}{\hat\Phi(q)}\left(1-\frac{\hat\psi'(0)\hat\Phi(q)}{q}\right)+
\frac{\hat\psi'(0)}{q}\frac{q-\hat\Phi(q)\psi'(\hat\Phi(q))}{q}.
\]

Now the claim follows from \eqref{eq:L}. 
\end{proof}

\subsection{Brownian motion}
\begin{thm}
\label{th.brownian}
Let $X$ be a standard Brownian motion $B\equiv\{B(t):t\in\rr\}$. Then,
for each $t>0$,
\[
\prob{M(t)>u}= \exp (-2u)
\left(
2(1+t)\Psi(\sqrt t)-\sqrt{\frac{2t}{\pi}}\exp\left( -\frac{t}{2} \right)
\right),
\]
where $\Psi(x)=\prob{\mathcal N>x}$ for a standard normal random variable $\mathcal N$.
\end{thm}
\begin{proof}
Since $B\in\clasS^-$, $Q(0)\de Q_e$ has an exponential distribution with mean $1/2$.
Thus,
\begin{align*}
\prob{M(t)>u}
&=\prob{Q(0)+\inf_{s\in[0,t]} (B(s)-s)>u}\\
&=
\int_u^\infty \prob{\inf_{s\in[0,t]} (B(s)-s)>u-x}2\exp (-2x)\D x\\
&=
2 \exp (-2u)\int_0^\infty \prob{\sup_{s\in[0,t]} (B(s)+s)<y} \exp (-2y)\D y\\
&=
\exp (-2u)
\ee \exp\left(-2\sup_{s\in[0,t]}(B(s)+s)\right)
\end{align*}
and the claim follows after some elementary computations (see also \cite[Eqn.\ (1.1.3)]{Borodin96} or
\cite{Baxter57}).
\end{proof}

\section{Asymptotics}
\label{sec:asymptotics}
In this section we consider the asymptotics of $\prob{M(T_u)>u}$
for a variety of functions $T_u$ and $u$ large.
As usual, heavy-tailed and light-tailed scenarios need to be
addressed separately.
\subsection{Heavy-tailed case}
In this section we shall work with the following assumption about the L\'evy process $X$:
\begin{assumption}
\label{ass:1}
For $\alpha>1$, let $X(1)\in \mathscr{RV}(-\alpha)$ --
the class of distributions with a complementary distribution function that
is regularly varying at $\infty$ with index $-\alpha$. Moreover, if $\alpha\in(1,2)$, then
in addition
\[
\lim_{x\toi}\frac{\prob{X<-x}}{\prob{X>x}}=\rho\in[0,\infty).
\]
\end{assumption}
We start with the following general proposition.
\begin{prop}
\label{prop:SLLN}
Let $X$ be a L\'evy process such that $\ee
X(1)=0$. Then, for any $\varepsilon>0$,
\[
\lim_{u\toi}\prob{\left|\frac{\underline{Y}(u)}{u}+1\right|>\varepsilon}=0.
\]
\end{prop}
\begin{proof}
First note that for such $X$ we have
\begin{equation}
\label{eq:SLLN}
\frac{X(t)}{t}\to 0 \text{ a.s. or equivalently } \frac{Y(t)}{t}\to -1 \text{ a.s.}.
\end{equation}
Fix $\varepsilon>0$, then
\[
\prob{\left|\frac{\underline{Y}(u)}{u}+1\right|>\varepsilon}
= \prob{\frac{1}{u}\inf_{t\in[0,u]} Y(t)>-1+\varepsilon}+\prob{\frac{1}{u}\inf_{t\in[0,u]} Y(t)<-1-\varepsilon}=:I_1(u)+I_2(u).
\]
Now \eqref{eq:SLLN} implies that
\[
I_1(u)\le \prob{\frac{Y(u)}{u}>-1+\varepsilon}\to 0.
\]
Observe that, for any $T\le u$,
\begin{align*}
I_2(u)&\le \prob{\frac{1}{u}\inf_{t\in[0,u]} Y(t)<-1-\varepsilon, \inf_{t\in[0,T]}Y(t)<\inf_{t\in[T,u]}Y(t)}\\
&\quad+
\prob{\frac{1}{u}\inf_{t\in[0,u]} Y(t)<-1-\varepsilon, \inf_{t\in[0,T]}Y(t)\ge\inf_{t\in[T,u]}Y(t)}\\
&\le 
\prob{\inf_{t\in[0,T]}Y(t)<\inf_{t\in[T,u]}Y(t)}
+
\prob{\frac{1}{u}\inf_{t\in[T,u]} Y(t)<-1-\varepsilon}=:I_{21}(u)+I_{22}(u).
\end{align*}
Again \eqref{eq:SLLN} implies that $I_{21}(u)\to 0$. As for $I_{22}(u)$, note that
\[
I_{22}(u)\le \prob{\frac{1}{u}\inf_{t\in[T,u]}X(t)<-\varepsilon}.
\]
Now we will show, that for any $\delta>0$ and $u$ large enough we have
\[
\prob{\frac{1}{u}\inf_{t\in[T,u]}X(t)<-\varepsilon}<\delta,
\]
which proves that $I_{22}(u)\to 0$ and completes the proof.

Indeed, by \eqref{eq:SLLN}, $T$ can be chosen such that
\[
\prob{\forall{t\ge T}:|X(t)|<\delta t}>1-\delta. 
\]
Therefore, for such $T$ and $\delta<\varepsilon$ we have
\begin{align*}
\prob{\frac{1}{u}\inf_{t\in[T,u]}X(t)<-\varepsilon}
&=
\prob{\frac{1}{u}\inf_{t\in[T,u]}X(t)<-\varepsilon, \forall{t\ge T}:|X(t)|\ge\delta t}\\
&\quad +
\prob{\frac{1}{u}\inf_{t\in[T,u]}X(t)<-\varepsilon, \forall{t\ge T}:|X(t)|<\delta t}\\
&< 
\delta,
\end{align*}
where we used the fact that the last probability equals zero.
\end{proof}

In the sequel we say that $f(n)\sim g(n)$ if $f(n)/g(n)\to 1$ as $n\toi$.
\begin{prop}
\label{prop.2}
Assume that  the L\'evy process $X$ satisfies \autoref{ass:1}.\\
(i) If $f(n)\ge n$, then
\[
\prob{X(n)>f(n)}\sim n\prob{X(1)>f(n)},
\] as $n\to\infty$, $n\in \nn$.\\
(ii) As $u\to\infty$,\[\prob{Q_e>u}\sim\frac{u}{\alpha-1}\prob{X(1)>u}.\]
\end{prop}
\begin{proof}
{\it Ad (i).} These asymptotics can be found in, e.g., \cite{Cline98} for $\alpha\ge 2$ and \cite{Borovkov00,Borovkov03} for $\alpha\in(1,2)$; see also \cite{Denisov08} for a recent treatment.
{\it Ad (ii).} See, e.g., \cite{Asmussen98,Dieker06}.
\end{proof}

We now state the main result of this subsection: the exact asymptotics of $\prob{M(T_u)>u}$.
\begin{thm}
\label{thm:heavy}
Assume that the L\'evy process $X$ satisfies \autoref{ass:1}. Then
\[
\prob{M({T_u})>u}\sim\prob{Q_e>u+{T_u}}+{T_u}\,\prob{X(1)>u+T_u},\as u.
\]
\end{thm}

The asymptotics in \autoref{thm:heavy} can be made more explicit. Part (ii) of \autoref{prop.2} immediately leads to the following corollary.
\begin{cor}
\label{cor:heavy}
Assume that the L\'evy process $X$ satisfies \autoref{ass:1}. Then
\[
\prob{M({T_u})>u}\sim
\left\{
\begin{array}{lcl}
\frac{1}{\alpha-1}u\,\prob{X(1)>u} & \text{when} & {T_u}=o(u),\\
\frac{A+\alpha}{\alpha-1}(A+1)^{-\alpha}{T_u}\,\prob{X(1)>{T_u}} & \text{when} & u\sim A{T_u},\\
\frac{\alpha}{\alpha-1}{T_u}\,\prob{X(1)>{T_u}} & \text{when} &
u=o({T_u}),
\end{array}
\right.
\]
as $u\toi.$
\end{cor}
\begin{proof}[{\bf Proof of \autoref{thm:heavy}}]
The proof consists of an upper bound and a lower bound. We use the notation $T_u^-:=\lfloor T_u\rfloor$ and $T_u^+:=\lceil T_u\rceil.$

\vb

{\it Upper bound}. To prove an (asymptotically) tight upper bound
for $\prob{M({T_u})>u}$, first we observe that for any
$\varepsilon>0$, using that $Q(0)\de Q_e$ is independent of $\{X(t):t\ge 0\}$,
\begin{align*}
\prob{M({T_u})>u} &\le \prob{M({{T_u^-}})>u}\:\le\:     \prob{Q_e+X({{T_u^-}})\ge u+ {{T_u^-}}}\\
&\le \prob{Q_e>(1-\varepsilon)(u+{{T_u^-}})}+\prob{X({{T_u^-}})>(1-\varepsilon)(u+{{T_u^-}})}\\
&\quad+\:\prob{Q_e>\varepsilon(u+{{T_u^-}})}\prob{X({{T_u^-}})>\varepsilon(u+{{T_u^-}})}\\
&=:\pi_{1}^+(u)+\pi_{2}^+(u)+\pi_{3}^+(u).
\end{align*}
Using (i) of \autoref{prop.2} and the strong law of large
numbers for $X$, it is easy to show that $\pi^+_{3}(u)=o(\pi_{1}^+(u))$
for a fixed $\varepsilon$. It is standard now to show that
\[
\lim_{\varepsilon\to0}\limsup_{u\toi}\frac{\pi_{1}^+(u)}{\prob{Q_e>u+T_u}}=1.
\]
Moreover,
\[
\lim_{\varepsilon\to0}\limsup_{u\toi}\frac{\pi_{2}^+(u)}{T_u\,\prob{X(1)>u+T_u}}=1,
\]
due to item (i) in \autoref{prop.2}. This establishes the upper
bound.

\vb

{\it Lower bound.} As for the lower bound observe that
\begin{align*}
\prob{M({T_u})>u}&\ge \prob{M({{T_u^+}} )>u}\:\ge\:
\prob{Q_e+\underline{Y}({{T_u^+}})>u,X({{T_u^+}})-{{T_u^+}}
               -\underline{Y}({{T_u^+}})<\varepsilon {{T_u^+}}}\\
&\ge\prob{Q_e+X({{T_u^+}})>u+(1+\varepsilon){{T_u^+}}}\prob{X({{T_u^+}})-{{T_u^+}}-\underline{Y}({{T_u^+}})
                 <\varepsilon {{T_u^+}}}\\
&=: \pi^-_{1}(u)\pi_{2}^-(u).
\end{align*}
By \autoref{prop:SLLN}, $\pi_{2}^-(u)\to 1$ as $u\toi$. Also,
\begin{align*}
\pi_{1}^-(u) &\ge
\prob{Q_e+X({{T_u^+}})-\varepsilon {{T_u^+}}/2>u+(1+\varepsilon/2){{T_u^+}},
X({{T_u^+}})>-\varepsilon {{T_u^+}}/2}\\
&\ge
\prob{\{Q_e>u+(1+\varepsilon/2){{T_u^+}},X({{T_u^+}})>-\varepsilon {{T_u^+}}/2\}\cup\{X({{T_u^+}})>u+(1+\varepsilon/2){{T_u^+}}\}}\\
&=
\prob{Q_e>u+(1+\varepsilon/2){{T_u^+}}}\prob{X({{T_u^+}})>-\varepsilon {{T_u^+}}/2}+\prob{X({{T_u^+}}))>u+(1+\varepsilon/2){{T_u^+}}}\\
&\quad
-\prob{Q_e>u+(1+\varepsilon/2){{T_u^+}}}\prob{X({{T_u^+}})>u+(1+\varepsilon/2){{T_u^+}}}\\
&=: \pi_{3}^-(u)\pi_{4}^-(u)+\pi_{5}^-(u)-\pi_{6}^-(u)\pi_{7}^-(u),
\end{align*}
where we again used that $Q(0)\de Q_e$ and $\{X(t):t\ge 0\}$ are
independent. By the strong law of large numbers, $\pi_{4}^-(u)\to1$ as
$u\toi$. Moreover, it is easy to show that
$\pi_{6}^-(u)\pi_{7}^-(u)=o(\pi_{3}^-(u))$. Now the lower bound follows by noting
that
\[
\lim_{\varepsilon\downarrow
0}\liminf_{u\toi}\frac{\pi_{3}^-(u)}{\prob{Q_e>u+{T_u}}}=1,
\]
and that (i) of \autoref{prop.2} yields
\[
\lim_{\varepsilon\downarrow 0}\liminf_{u\toi}\frac{\pi_{5}^-(u)}{{T_u}\,\prob{X(1)>u+{T_u}}}=1.
\]
This completes the proof.
\end{proof}

\subsubsection{Stable L\'evy processes}
Following the notation from \cite{Samorodnitsky94}, let $S_{\alpha}(\sigma,\beta,\mu)$
be a stable law with index $\alpha\in(0,2)$, scale parameter $\sigma>0$, skewness parameter $\beta\in[-1,1]$
and drift $\mu\in\rr$. We call $X$ an $(\alpha,\beta)$-stable L\'evy process if
$X$ is a Levy process and
$X(1)$ has the same distribution as $S_\alpha(1,\beta,0)$.

Let
\[
B(\alpha,\beta):=\frac{\Gamma(1+\alpha)}{\pi}\sqrt{1+\beta^2\tan^2\left(\frac{\pi\alpha}{2}\right)}
\sin\left(\frac{\pi\alpha}{2}+\arctan\left(\beta\tan\left(\frac{\pi\alpha}{2}\right)\right)\right),
\]
and let $X$ be an $(\alpha,\beta)$-stable L\'evy process with $\alpha\in(1,2)$ and $\beta\in (-1,1]$.
Then,
\[
\prob{X(1)>u}\sim \frac{B(\alpha,\beta)}{\alpha}u^{-\alpha},
\]
see, e.g., \cite[Prop.\ 2.1]{Port89}.
Now \autoref{thm:heavy} can be rephrased as follows.
\begin{cor}
For an $(\alpha,\beta)$-stable L\'evy process $X$ with $\alpha\in(1,2)$ and $\beta\in (-1,1]$,
\[
\prob{M({T_u})>u}\sim
\left\{
\begin{array}{lcl}
\frac{1}{\alpha-1}\frac{B(\alpha,\beta)}{\alpha}u^{1-\alpha} & \text{when} & {T_u}=o(u),\\
\frac{A+\alpha}{\alpha-1}(A+1)^{-\alpha} \frac{B(\alpha,\beta)}{\alpha}{T_u}^{1-\alpha} & \text{when} & u\sim A{T_u},\\
\frac{\alpha}{\alpha-1} \frac{B(\alpha,\beta)}{\alpha}{T_u}^{1-\alpha} & \text{when} & u=o({T_u}),
\end{array}
\right.
\]
as $u\toi.$
\end{cor}
\subsection{Light-tailed case}
In this subsection, we consider the light-tailed situation, also frequently referred to as the Cram\'er case. Throughout, with $\phi({\vartheta}):=\log\ee\exp({\vartheta} X(1))$ denoting the {\it cumulant function}, we impose the following assumption.
\begin{assumption}
\label{assum:2}
Let
\[
\beta^{\star}:=\sup\{\beta:\ee e^{\beta X(1)}<\infty\}
\]
Assume that $\beta^{\star}>0$ and there exists
${\vartheta}^{\star}\in(0,\beta^{\star})$ such that
$\phi({\vartheta}^{\star})={\vartheta}^{\star}$. 
\end{assumption}
For $r\ge 0$,
define
\[
I(r):=\sup_{{\vartheta}>0} \left({\vartheta} r-\phi({\vartheta})\right).
\]
\begin{prop}
\label{Cramer}
Under \autoref{assum:2}, the following statements hold.

(i) As $u\toi$,
\[
\log\prob{Q_e>u} \sim-{\vartheta}^{\star}u.\]
(ii) For all $u>0$,
\[
\prob{Q_e>u} \le e^{-{\vartheta}^\star u}.
\]
(iii) The function $I$ obeys
\[
I(1)<\infty \quad \text{and}\quad I'(1)\le{\vartheta}^\star.\]
(iv) For any $\varepsilon>0$,
\[
\liminf_{u\toi}\frac{1}{u}\log\prob{\underline{Y}(u)>-\varepsilon u}\ge - I(1)\]
\end{prop}
\begin{proof} For (i) and (ii), we refer to \cite[Thms. 5.1 and 5.2 of Chapter XIII]{Asmussen02} or \cite[Proposition 5.1 and Remark 5.1]{Debicki10}.
For (iii), notice that
$I(1)=\sup_{{\vartheta}>0}({\vartheta}-\psi({\vartheta}))$ is attained for
${\vartheta}\in(0,{\vartheta}^\star)$; therefore also $I'(1)\le{\vartheta}^\star$.
As for (iv), observe that
\[
\prob{\underline{Y}(u)>-\varepsilon u}
=
\prob{\frac{Y(u\,\cdot)}{u}\in A_{\varepsilon}},
\]
where
\[
A_\varepsilon:=\left\{f\in D[0,1]: f(t)>-\varepsilon, \forall{t\in [0,1]}\right\}
\]
and $D[0,1]$ is the space of c\`adl\`ag functions on $[0,1]$.
Using sample-path large deviations results for L\'evy processes
(see \cite[Theorems 5.1 and 5.2]{DeAcosta1994}),  we now obtain that
\[
\liminf_{u\toi}\frac{1}{u}\log\prob{\underline{Y}(u)>-\varepsilon u}
\ge -\inf\{\psi(f):f\in A_\varepsilon\cap C[0,1]\},
\]
where $\psi(f):=\int_0^1 I(f'(t)+1)\,\D t.$
Now observe that the path $f^\star\equiv 0$ is in $A_\varepsilon.$
The stated follows by realizing that $\psi(f^\star)= I(1).$\end{proof}

Now we can proceed with the main result of this subsection.
\begin{thm}
\label{thm:lt}
Assume that the L\'evy process $X$ satisfies \autoref{assum:2}. Then
\[
\log \prob{M({T_u})>u} \sim -u {\vartheta}^\star - {T_u}I(1),\as u.
\]
\end{thm}

The asymptotics in \autoref{thm:lt} can trivially be made more explicit by comparing both exponential decay rates.
The intuition behind the following corollary is that, in large deviations language, the most likely path corresponding to the rare event under study first builds up from an empty system to level $u$ (at time 0), and then remains at level $u$ for the nest $T_u$ time units; both parts of the path result in both contributions to the decay rate (i.e., $-u\vartheta^\star$ and $-T_uI(1)$). Then, depending on whether $T_u$ is small or large with respect to $u$, one of these two contributions dominates.
\begin{cor}
\label{cor:lt}
Assume that the L\'evy process $X$ satisfies \autoref{assum:2}. Then
\[
\prob{M({T_u})>u}\sim
\left\{
\begin{array}{lcl}
-u {\vartheta}^\star  & \text{when} & {T_u}=o(u),\\
-T_u(A\vartheta^\star +I(1))& \text{when} & u\sim A{T_u},\\
-{T_u}I(1) & \text{when} &
u=o({T_u}),
\end{array}
\right.
\]
as $u\toi.$
\end{cor}
\begin{proof}[{\bf Proof of \autoref{thm:lt}}]
The proof again consists of two bounds.

\vb

{\it Lower bound.} Observe that the probability of interest is,
for any $\varepsilon>0$, bounded from below by
\[
\prob{Q_e>u+\varepsilon {T_u}}\prob{\underline{Y}({T_u})> -\varepsilon {T_u}}.
\]
Now the lower bound follows by combining parts (i) and (iv) of \autoref{Cramer}, and
then sending $\varepsilon\downarrow 0$.

\vb

\newcommand{\IT}[1]{{\mathbb I}(#1)}

{\it Upper bound.}
Observe that
\[
\underline{Y}(t)\le \left(X(t)-{t}\right) \IT t,
\]
where $\IT t$ denotes the indicator function
$1_{\{X(t)/{t}\in(0,1)\}}.$ Thus,
\begin{align*}
\prob{M({T_u})>u} &\le \prob{Q_e+(X({T_u})-{T_u})\IT{T_u} >u }\\
&=\int_\rr\prob{Q_e>u-x{T_u}+\IT{T_u}}\D\prob{\frac{X({T_u})}{{T_u}}\IT{T_u}\le x}\\
&=\int_0^1\prob{Q_e>u-x{T_u}+{T_u}}\D\prob{\frac{X({T_u})}{{T_u}}\IT{T_u}\le x}\\
&\le e^{-{\vartheta}^\star u}\int_0^1 e^{-{\vartheta}^\star
{T_u}(1-x)}\D\prob{\frac{X({T_u})}{{T_u}}\IT{T_u}\le x},
\end{align*}
where the last inequality follows from part (ii) of \autoref{Cramer}. The
sequence $\{X(u)\IT{u}/u\}$ satisfies the large deviations principle on $((0,1),\mathcal B(0,1))$
with rate $u$ and rate function $I(\cdot)$. Thus, Varadhan's lemma
\cite[Theorem.\ 4.3.1]{Dembo98} implies
\[
\lim_{u\toi}\frac{1}{{T_u}}\log\int_0^1 e^{-{\vartheta}^\star {T_u}(1-x)}\D\prob{\frac{X({T_u})}{{T_u}}\IT{T_u}\le x}
=
-\inf_{x\in(0,1)}\left({\vartheta}^\star(1-x)+I(x)\right)=I(1),
\]
where the last equality is due to part (iii) of \autoref{Cramer}
and convexity of $I(\cdot)$.
\end{proof}

\twocolumn
\small
\bibliography{Infimum}
\end{document}